\documentclass[hyperref]{lmcs}
\pdfoutput=1

\usepackage{lastpage}
\renewcommand{\dOi}{14(3:5)2018}
\lmcsheading{}{1--\pageref{LastPage}}{}{}%
{Dec.~20,~2017}{Jul.~31,~2018}{}

\usepackage{amssymb}
 \usepackage{amsthm}
 \usepackage{latexsym}
\usepackage{amsmath}

\usepackage{csquotes}
\def\signed #1{{\leavevmode\unskip\nobreak\hfil\penalty50\hskip2em
  \hbox{}\nobreak\hfil(#1)%
  \parfillskip=0pt \finalhyphendemerits=0 \endgraf}}

\newsavebox\mybox
\newenvironment{aquote}[1]
  {\savebox\mybox{\emph{#1}}\begin{quote}}
  {\signed{\usebox\mybox}\end{quote}}

\numberwithin{equation}{section}

\newcommand{\Iff}{\mbox{$\Longleftrightarrow$}}

\def\proof{{\bf Proof.}\ }

\def\ie{\emph{i.e.}}
\def\pes{\emph{e.g.}}
\def\Pes{\emph{E.g.}}

\def\F{{\mathcal F}}

\def\P{{\mathcal P}}

\def\U{{\mathcal U}}

\def\V{{\mathcal V}}

\def\C{{\mathcal C}}

\def\N{{\mathbb N}}

\def\R{{\,\mathbb R}}

\def\hR{{}^{*}\R }

\def\hM{{}^{*}\! M }
\def\hP{{}^{*}\! P }
\def\hRR{{}^{*}\! R }

\def\g{{}^{*}\! g }

 \def\p{{}^{*}\! p}
 \def\q{{}^{*}\! q}
 \def\f{{}^{*}\! f}
 
 \def\X{{}^{*}\! X}
 \def\hA{{}^{*}\! A}
 \def\ns{{}^{*}\!}
 \def\hB{{}^{*}\! B}

	\def\B{\mathcal{B}}

\def\bg{{\beta}}

\def\yg{{\upsilon}}

\def\+#1{\vec{#1}}

\def\uf{{\mathfrak{u}}}

\def\*{\times}
\def\st{{}^{\star}\!\, }
\def\0{\emptyset}
\def\7{\setminus}
\def\_{\overline}

\def\<{\prec}

\def\incl{\subseteq}

\def\transfer{\textit{Transfer Principle}}
\def\ind{\textsf{Ind}}
\def\poss{\textsf{Poss}}
\def\tran{\textsf{Tran}}
\def\zfc{\textsf{ZFC}}

\def\qed{${}$\hfill $\Box$}

\def\trans{transfer principle}

\def\Cl{{\C\ell}}
\def\transfer{\textit{Transfer Principle}}
\def\equ{\textsf{(equ)}}
\def\dir{\textsf{(dir)}}
\def\comp{\textsf{(comp)}}
\def\ind{\textsf{Ind}}
\def\poss{\textsf{Poss}}
\def\tran{\textsf{Tran}}

\title[Leibniz's Principles  within Functional Extensions]{A topological interpretation of three Leibnizian Principles 
within the Functional Extensions}
\thanks{It is a great pleasure to  dedicate this paper to Furio: how many joint papers in the Eighties and Nineties deal with mathematical structures arising from fundamental philosophical principles!}

\author[Marco Forti]{Marco Forti}
\address{Dipartimento di Matematica\\
Universit\`{a} di Pisa\\
Via Buonarroti 1C\\
56100 Pisa}

\email{marco.forti@unipi.it}
\keywords{transfer principle, indiscernibles, nonstandard models, 
functional extensions}

\copyrightinfo{}

\begin{document}

\begin{abstract}
    
       Three philosophical principles are often quoted in connection with Leibniz:
     ``objects sharing the same properties are the same object'' (\emph{Identity of indiscernibles}), 
     ``everything can        
       possibly exist, unless it yields contradiction'' (\emph{Possibility as consistency}), and
     ``the ideal elements          
      correctly determine the real things'' (\emph{Transfer}).

      Here we give a precise logico-mathematical formulation 
       of these principles within the framework of
      the Functional Extensions, mathematical structures     
       that generalize at once compactifications, completions, and 
       elementary          %
       extensions of models. In this context, the above Leibnizian principles 
	appear  
       as topological or algebraic  properties, namely: a property of separation, 
       a property of compactness, and 
       a property of directeness, respectively. 
      
   Abiding by this interpretation, we obtain the somehow surprising 
   conclusion that these Leibnizian principles 
    \emph{may be fulfilled in pairs, but not all          
       three together}. 
\end{abstract}

\maketitle

%\subjclass{03A05, 03H05, 03E65, 54D35, 54D80}

%% be sure to place the \maketitle command here.

%\newpage

\section*{Introduction}

The great philosopher, mathematician and logician Gottfried Wilhelm von Leibniz, in many of his philosophical writings, has been inspired by, and consequently has given inspiration to several important mathematical ideas. In this article we consider, within a mathematical framework,
three philosophical principles that are often quoted in connection with Leibniz.\footnote{For the Leibniz's quotes see \href{http://www.leibnizedition.de/}{www.leibnizedition.de/}. }
We cannot discuss and articulate the philosophical aspects of these principles: we simply give here a few quotes in order to justify and explain our mathematical interpretation of these Leibnizian principles.
\smallskip

\textbf{Identity of indiscernibles}
      \begin{center}
          \emph{``objects sharing the same properties are the same object''}
      \end{center} 
 
     \begin{aquote}{Monadology}
    
    There are never in nature two beings which are perfectly identical to 
    each other, and in which it is impossibile to find any internal difference
     \ldots
   \end{aquote}
   
   \begin{aquote}{Fourth letter to Clarke}
     Two  indiscernible individuals cannot exist. [\ldots]
     To put two 
     indiscernible things is to put the same thing under two names.
   \end{aquote}
   \medskip
   \begin{aquote}{Fifth letter to Clarke}
     \textelp{} dans les choses sensibles on n'en trouve 
     jamais deux indiscernables \textelp{}
   \end{aquote}
\medskip

\textbf{Possibility as consistency}
\begin{center}
          \emph{``everything can        
        possibly exist, unless it yields contradiction''}
     \end{center}
     
     \begin{aquote}{\emph{unpublished}, 1680 ca.}
       \begin{itemize}[-]
         \item \emph{Impossible} is what yields an absurdity.
         \item \emph{Possible} is not impossible.
         \item \emph{Necessary} is that, whose opposite is impossible.
         \item \emph{Contingent} is what  is not  necessary.
         \end{itemize}
       \end{aquote}   
   \medskip
   \begin{aquote}{Dialogue between Theofile and Polydore}
    \textelp{} nothing is absolutely necessary, when the contrary 
     is possible. \textelp{}\\
    Absolutely necessary is \textelp{} that whose opposite yields a contradiction. 
  \end{aquote}
 
 \medskip
    
     \textbf{Transfer principle}
    \begin{center}
         \emph{``the ideal elements          
           correctly determine the real things''}
     \end{center}
     
      \smallskip
 \begin{aquote}{\emph{Letter to Johann Bernoulli}, 1698}
     Perhaps \emph{the infinite and infinitely small} [numbers] 
    that we conceive \emph{are imaginary, nevertheless} [they are] \emph{suitable 
    to determine the real things}, as usually do the imaginary roots.
    They are situated in the ideal regions, from  where things are ruled by
    laws, even though they do not lie in the part of matter.
  \end{aquote}

\bigskip  
  
  In this paper, we try and give a precise mathematical formulation 
  of these principles in the context of
  the \emph{Functional Extensions} of \cite{F12}, structures     
   which generalize at once compactifications and completions of topological spaces, and nonstandard    extensions of      %
   models (see also \cite{BDF06}). 
    Given a set $M$, a 
   \emph{functional extension} of $M$ is a superset $\hM$ of
   $M$ such that every function $f:M\to M$ has a distinguished 
   extension $\f: \hM\to\hM$ that preserves compositions and equalizers. 
   Moreover, assuming that the so called Puritz preorder of $\hM$ is directed, 
   the operator $\ast$ can be appropriately defined so as to 
   provide also
   all properties $P$ and all relations  $R$ on $M$ with  distinguished 
   extensions $\hP$ and $\hRR$ on $\hM$.
   
  Following the basic idea that the elements of the [``standard''] set $M$ are 
  the ``real objects'' 
  of the ``actual world'', whereas the [``nonstandard"] extension $\hM$ 
  contains also the ``ideal elements'' of all ``possible worlds'', 
  an appropriate interpretation of 
  the Leibniz's principles   in the context of functional extensions
  might be the following:
  
  \begin{description}
     \item[Ind] \emph{different elements of $\,\hM$ are \emph{separated} by the 
     extension $\hP$ of some ``real" property $P$ of $M$;} 
 
     \item[Poss] \emph{a family of ``real" properties of $M$ that are not \emph{contradictory} in $M$ has  extensions to $\hM$ that
      are  \emph{all simultaneously 
     satisfied} in $\hM$};  
     \item[Tran] \emph{a statement involving ``real" elements, properties 
     and relations of $M$ is \emph{true in $M$ }if and only if the corresponding 
     statement about their ``ideal" extensions is \emph{true in 
     $\hM$}.}
     (Clearly here one has to admit only ``first order''
     statements, so as to avoid trivial inconsistencies.)
 \end{description}
  
  We shall see below that the ``real" properties of $M$, when extended to $\hM$, may be taken as the 
  \emph{clopen} subsets of a topology  on $\hM$ (the so called $S$-\emph{topology}), and
  so the above principles \ind\ and \poss\ turn out to be respectively 
    a
   property of \emph{separation} and a property  of \emph{compactness}
   of the $S$-topology of $\hM$. On the other hand, the principle \tran\ requires \emph{all} the conditions postulated for the $\ast$-extension in  \cite{F12}, including the order-theoretic property of directeness, so as to take care of properties and relations of arbitrary arities.
 
 On the basis  of results similar to those of \cite{DF05,BDF06}, we obtain 
the somehow 
surprising consequence that these Leibnizian principles 
\emph{can be fulfilled in pairs, but not all          
   three together}. 
  Grounding on the mathematical results of this paper, and, chiefly, on the philosophical writings of Leibniz himself,  \pes
 
 \begin{quote} 
 Les parties du temps ou du lieu [...] sont des choses ideales, ainsi
elles se rassemblent parfaitement comme deux unit\'es abstraites. 
Mais il n'est pas de m\^eme de deux Uns concrets [...] c'est \`a 
dire veritablement actuels.\footnote{~The 
parts of time or place \ldots are ideal things, so 
they perfectly resemble like
two abstract unities. But it is not so with two concrete Ones, \ldots
that is truly actual [things].
} (\emph{Fifth letter to Clarke})

\end{quote} 
 we 
  decide to call \emph{Leibnizian} those functional extensions that satisfy \poss\ together with \tran, letting hold the principle \ind\ only inside the ``standard" universe.
  
\medskip
The paper is organized as follows:
\begin{itemize}
\item In Section \ref{top}, we recall the precise definition of 
\emph{functional
extension} together with the main  properties stated in \cite{F12}.  
In particular, we establish the fundamental result that the functional extensions are \emph{all and only the complete elementary} (=nonstandard) extensions of any model $M$.
\item In Section \ref{stop} we introduce the $S$-topology on functional extensions, and we determine which functional
extensions satisfy the \trans\  \tran\ together with either \poss, or \ind, according to the  properties of the $S$-topology.
Consequently, in Subsection \ref{tros} we see  how  to obtain  
functional
extensions where the principles \poss\ and \tran\ hold {together}, whereas in Subsection \ref{trin} we examine the structure of the extensions verifying \ind\ and \tran\ together. The 
\emph{impossibility of satisfying simultaneously the three Leibnizian
principles}
follows immediately.
\item In Section \ref{shaus}
we consider which set theoretic hypotheses are needed  in order to get functional extensions satisfying the principle \ind. Leaving apart these conditions, in Subsection \ref{lchoi} we motivate by his own quotes the (supposed) preference of Leibniz himself for dropping off \ind\ while maintaining the remaining principles \poss\ and \tran. 
\item In Section \ref{inos} we suggest how a slight weakening of the notion of functional extension allows for admitting the simultaneous holding of the principles \ind\ and \poss, still preserving a large amount of \tran, (by dropping off only the incompatible``analytic" part of the condition \equ\ that forces everywhere different functions on $X$ to have everywhere different extensions to $\X$).
\item A few concluding remarks and three connected set theoretic open questions  can be found in the final 
Section \ref{froq}.
\end{itemize}

\smallskip
In general, we refer to \cite{Eng} for all the topological notions and
facts used in this paper, and to \cite{CK} for definitions and facts
concerning ultrapowers, ultrafilters, and nonstandard models. General 
references for nonstandard Analysis could be \cite{Ke76,nsa}; 
specific for our ``elementary'' approach is  \cite{BDF06}.

\section{Functional Extensions and  the \trans}
\label{top}
The Transfer Principle  \tran\ is the very ground of the usefulness of the 
nonstandard methods in mathematics. It allows for obtaining correct 
results about, say, the real numbers by using \emph{ideal elements}, like 
actual \emph{infinitesimal} or \emph{infinite} numbers.
In this section we review the main features of the \emph{functional extensions}  
introduced in the paper \cite{F12} (see also \cite{BDF06}), with the goal of characterizing all  \emph{nonstandard } (= complete elementary) extensions by means of a few simple properties of an operation $\ast$ that assigns an appropriate extension $\f:\X\to\X$ to  each function $f:X \to X$.

These structures are a sort of 
%Stone-$\check{\rm C}$ech  
\emph{compactifications or completions of 
discrete spaces}  
that  turn out to comprehend \emph{all and only} the \emph{nonstandard extensions of models}.

Recall that the \emph{Puritz (pre)ordering} $\lesssim_P$  of $\X$ is defined by $\eta\lesssim_P \xi$ if there exists $f:X\to X$ such that $\f(\xi)=\eta$ (see \cite{pu,ng}).
 Then we set our fundamental definition as follows:
\begin{defi}\label{text}
A
\emph{functional extension} of $X$ is  a  proper superset $\X$ of $X$ where a
\emph{distinguished extension} $\f:\X \to \X$ is assigned to
each function $f:X \to X$, so as to satisfy the following conditions
for all $f,g:X \to X$. 
%and all $\xi\in\X$:}
\begin{description}
\item[$(\mathsf{comp})$]  \textbf{Preservation of compositions:}\ \  $ \ns(g\circ f)= \g\circ\f$\\
\ \ \ie\ \ \ $\g( \f(\xi)) = \ns(g\circ f) (\xi)$ \ for all\  $\xi\in\X;\ $
%\bigskip
\item[$(\mathsf{equ})$] \textbf{Preservation of equalizers:}\footnote{~
In \cite{F12}, this property is denoted by $\mathsf{(diag)}$, and called \emph{preservation of the diagonal}, in view of the fact that $\chi_{fg}(x)=\chi_{\Delta}(f(x),g(x))$, where $\chi_{\Delta}$ is the characteristic function of the diagonal $\Delta\incl X\times X$.
(For convenience we always assume that $0,1\in X$.)}
~$\ns(\chi_{fg})= \chi_{\f\g}$

		     where $\, \chi_{\phi\psi}$ is the characteristic
function of the equalizer $Eq(\phi,\psi)$\\ of the functions $\phi$ and $\psi$, \ie \ 
$\chi_{\phi\psi}(z)=\begin{cases}
     1 &if\ \ \phi(z)=\psi(z) , \\
      0 & {otherwise};
\end{cases}$
\end{description}
%Recall that the \emph{Puritz (pre)ordering} $\lesssim_P$  of $\X$ is defined by $\eta\lesssim_P \xi$ if there exists $f:X\to X$ such that $\f(\xi)=\eta$ (see \cite{pu,ng}).
%So we postulate
\begin{description}
\item [$(\mathsf{dir})$] \textbf{Directness of the Puritz order:}~
for all\  $\xi,\eta \in \X$\  there
exist  \\ $p_1,p_2:X\to X$ and $\zeta \in \X$ such that
$ \xi = \p_1(\zeta)$\ \ and\  \ $\eta = \p_2 (\zeta)$, \\ so that\ \ $\xi,\eta\lesssim_P \zeta$.
\end{description}
\end{defi}

The importance of the  property $\mathsf{(dir)}$, called
\textit{coherence} in \cite{DF05}, is due to the fact that, by providing an ``internal coding of
pairs", it allows
for extending  multivariate functions ``parametrically'': this
possibility is essential in order to get the full principle  \tran, which 
involves properties, relations, and functions of any arities. 
More precisely, the following facts that hold in every   functional
extension $\X$ of $X$ allow for considering only \emph{unary} functions
(see Subsection 3.2 of \cite{F12}):\\
\emph{- For all $\xi_{1},\ldots,\xi_{n}\in \X$ there exist
$ p_{1},\ldots,p_{n}: X \to X$ and  $\zeta\in \X$ such that
$\p_{i}(\zeta) = \xi_{i}$.}\\
\emph{- If $p_{1},\ldots,p_{n},q_{1},\ldots, q_{n}: X \to X$ and
$\xi,\eta \in \X$ satisfy  $\p_{i}(\xi) = \q_{i}(\eta)$, then\\ 
${}~~~~\ns (F \circ (p_{1},\ldots,p_{n}))(\xi) =
\ns (F \circ (q_{1},\ldots, q_{n}))(\eta)$ for all $F: X^{n} \to X$}.

\noindent It follows that there is a unique way of assigning an
extension $\ns F$ to every function $F:X^{n} \to X$ in such
a way that all compositions are preserved. Then,
by using the characteristic functions in $n$ variables, one can assign
an extension $\ns R$ also to all $n$-ary relations $R$ on $X$.

Several
important cases of the
\trans\ are easy to deduce:
\pes, if $f$ is constant, or injective, or surjective, or characteristic, then
so is  $\f$. In particular
 the  extension $\ns\chi_A$ of the characteristic function $\chi_A$ of a subset $A\incl X$, can be taken as the characteristic function of
the $\ast$-{\em extension} $\hA$ of
$A$ in $\X$, thus obtaining   a 
\emph{Boolean isomorphism} between the field $\P(X)$ of all subsets of $X$ 
and a field $\C\ell(\X)$ of subsets of $\X$ (see Theorem 1.2 of \cite{F12}).

      While both properties $(\mathsf{comp})$ and $(\mathsf{equ})$ of 
Definition \ref{text} are clear
instances of the \trans,  as they correspond  to the 
statements
$$\forall x\in X\, .\, f(g(x))=(f\circ g)(x)\ \
\mbox{and}\ \
\forall x\in X\, .\, f(x)=g(x)\ \Iff\ \chi_{fg}(x)=1, $$
respectively, the same is not apparent for the third property  $(\mathsf{dir})$, which is given in a \emph{second order} formulation. 
On the contrary,  an even stronger 
\emph{uniform version} of $(\mathsf{dir})$ can be obtained by \tran: simply take $p_1,p_2$ to be the compositions of a fixed bijection  $\delta :X\to X\times X \ $ 
with the ordinary projections $\pi_{1},\pi_{2}: X\times X \to X$,
and  apply \tran\ to the statement
$$\forall x,y\in X\, .\,  \exists z
\in X\, .\ p_1(z) =x,\ p_2(z) =y.$$

Hence, when $\X$ is a nonstandard extension of $X$, the three defining properties 
$(\mathsf{comp})$, $(\mathsf{equ})$, and $(\mathsf{dir})$ of  functional extensions 
are fulfilled by hypothesis. 
Conversely, we devoted \cite{F12} to prove  the 
  fact (partly anticipated in \cite{BDF06}) that the combination of \emph{three natural, simple instances} 
 of the \trans, namely $(\mathsf{comp}),(\mathsf{equ})$, and $(\mathsf{dir})$,
makes any functional  extension $\X$ a \emph{limit ultrapower} of $X$, thus 
providing the full \transfer\ \tran. So the functional extensions are
exactly the \emph{hyper-extensions} in the sense of 
\cite{BDF06} (\ie\ all nonstandard models).

We state without proof the corresponding theorem, and we direct the reader  to \cite{F12}, where a \emph{purely algebraic} proof  is given in all details, and two more (a \emph{topological} and a\emph{ purely logical} one) are outlined.

\begin{thm}[Thm.~2.2 and Cor.~2.3 of \cite{F12}]\label{eltop}~ 
Any nonstandard extension $\,\X$ of $X$ is a functional extension of $X$, and conversely
any functional extension $\X$ of $X$ satisfies the \trans\ \emph{\tran}.
\qed
\end{thm}

\medskip
\section{The \texorpdfstring{$S$}{S}-topology and the principles \ind\ and \poss.}
\label{stop}
In order to study our formalizations of the Leibnizian principles \ind\ and \poss\ within the
 functional 
extensions, it is natural to consider on $\X$ a topology that corresponds to the 
 classical \emph{$S$-topology} of  Nonstandard
 Analysis, already considered since \cite{rob66} 
 ($S$ stands for \emph{Standard}). In this topology, the \emph{closure} in $\X$ of a subset $A\incl X$ is given by its extension $\hA$, so the field of subsets $\Cl(\X)=\{\hA\incl\X\mid A\incl X\}$ is both  an \emph{open basis } of the $S$-topology and  the \emph{field of all its clopen} (=closed and open) subsets (since $\ns(X\7 A)=\X\7\hA$).
 
 Remark that   all functions 
$\f$ are \emph{continuous}  with respect to the $S$-topology, 
because $\f^{-1}(\hA)=\ns(f^{-1}(A))$ for all $A\incl X$, so the inverse images of clopen sets are clopen.
Hence for $\X$ being a \emph{topological extension} in the sense of \cite{DF05} it needs only that the $S$-topology be $T_1$ (and so, according to the next theorem, Hausdorff).

We begin by characterizing the \emph{separation properties} of the $S$-topology
in the same way as in Theorem 1.4 of \cite{DF05}:

\begin{thm}\label{sep}
    Let $\X$ be a functional extension of $X$, and put on $\X$ the $S$-topology. Then, as a topological space,
%\begin{\begin{itemize}
  $\, \X$ is $0$-dimensional, and either totally disconnected or
not $T_{0}$.
%\end{itemize}}
 It follows that  
 $\X$ with the $S$-topology is Hausdorff if and only if it is $T_{0}$.
\end{thm}

\proof
The $S$-topology has a clopen basis by definition. In
this topology the closure of the point $\, \xi$ is $\, M_{\xi} = \bigcap
_{\xi \in A} \hA$, an intersection of clopen sets.
If  $M_{\xi} = \{\xi \}$ for all $\,\xi\in\X$, then different points have disjoint clopen neighborhoods, hence the $S$-topology is regular (and totally disconnected).

Assume instead that there exist points  $\eta \ne \xi$ such that  $\eta\in M_{\xi}$.
Then the
$S$-topology is not $T_{0}$, because $\eta$ belongs to the same basic open (clopen) sets as $\xi$. In fact, given
$A\incl X$, $\xi \in \hA $ implies
         $\eta \in \hA $, by the
choice of $\eta$. Similarly  $\xi \notin \hA $ implies
$\xi \in\! \ns(X\7 A) $, hence  $\eta \in\! \ns(X\7 A) $ and so
$\eta \notin \hA $.

 In particular the $S$-topology is Hausdorff (in fact regular and 
totally disconnected)
    whenever it is $T_{0}$. 
    \qed

\bigskip
Now the principle \ind\  simply means that, given $\xi\ne\eta\in\X$ there exist disjoint subsets $A,B\incl X$ such that $\xi\in\hA$ and $\eta\in\hB$, \ie\
that  the $S$-topology of $\X$ is Hausdorff. 

On the other hand, the principle 
\poss\ states that all the $\ast$-extensions of a family of finitely compatible properties of $X$ are simultaneously satisfied in $\X$.  So the corresponding family $\F$ of subsets of $X$  has the finite intersection property, hence the intersection $\bigcap_{A\in\F} \hA$  of the corresponding $\ast$-extensions has to be nonempty. But requiring this is equivalent to require that every \emph{proper filter of clopen sets} in $\X$ 
has \emph{nonempty intersection}, \ie\ that the $S$-topology of $\X$ is 
\emph{quasi-compact}. (Following \cite{Eng}, we call \emph{compact} only 
Hausdorff spaces.)

So we can completely determine, from the $S$-topology, when the principles \ind\ and \poss\ hold in a functional extension:

\begin{cor}\label{ip}
 Let $\X$ be a functional extension of $X$ with the $S$-topology. Then
 \begin{enumerate}
     \item  the principle \emph{\textbf{Ind}} holds if and only if $\,\X$ is 
     Hausdorff;
  \item  the principle \emph{\textbf{Poss}} holds if and only if  $\,\X$ is 
quasi-compact; hence
  \item  the principles \emph{\textbf{Ind}} and  \emph{\textbf{Poss}} 
     hold simultaneously in 
     $\X$ if and only if 
     it is compact.  \qed
       \end{enumerate}
       
        \end{cor}
        
\subsection{Combining \tran\ with \poss.}\label{tros}
In order to combine the principle \tran\ with \poss, we recall that a nonstandard model whose $S$-topology is 
quasi-compact is commonly called an \emph{enlargement}. It is well known 
that every structure has \emph{arbitrarily saturated} 
 elementary extensions (see \pes\ \cite{CK}), and  any 
 $(2^{|X|})^{+}$-\emph{saturated elementary extension} of $X$ is easily seen to be a \emph{nonstandard enlargement} (see \pes\
 \cite{nsa} or \cite{BDF06}). Therefore Corollary \ref{ip} yields

\begin{thm}\label{start}
The nonstandard enlargements of $X$ are exactly those functional extension  of $X$ 
that satisfy both principles \emph{\tran}
and \emph{\poss}. \qed
\end{thm}

Thus 
 we get a lot 
 of functional extensions satisfying simultaneously \tran\ and \poss.
 
 \subsection{Combining \tran\ with \ind.}\label{trin}
 We pass now  to characterize all those functional
extensions of $X$ that satisfy simultaneously the principles \tran\ and \ind. 

The $S$-topology of these extensions is Hausdorff, 
according to Corollary \ref{ip}, so we know that the intersection $M_\xi=\bigcap_{\xi\in A}\hA$ is the singleton $\{\xi\}$, hence the filter $\U_\xi$ of the subsets $A$ of $X$ such that $\hA$ contains $\xi$ uniquely determines the point  $\xi\in\X$. Moreover, the filter  $\U_\xi$ is actually an \emph{ultrafilter} on $X$ such that
$$\U_{\f(\xi)}=\_f(\U_\xi)=\{ A\incl X\mid f^{-1}(A)\in\U_\xi \}\ \ for\ all\ f:X\to X.$$

So, when $\X$ is Hausdorff, we must have that $$\_f(\U_\xi)=\_g(\U_\xi)\ \Iff\ 
\f(\xi)=\g(\xi)$$ \ for\ all\ $\xi\in\X$ and all $ f,g:X\to X.$

Now $\f(\xi)=\g(\xi)$ holds if and only if $\chi_{\f\g}(\xi)=\ns(\chi_{fg})(\xi)=1$, or equivalently if and only if the equalizer $E(f,g)$  of $f$ and $g$ belongs to $\U_\xi$.

\smallskip
Recall that an ultrafilter $\U$ on $X$ is called \emph{Hausdorff}\footnote{~The
      property $(\mathsf{H})$ has been introduced in \cite{dt} under
      the name $(C)$. Hausdorff ultrafilters are studied in \cite{kt, DF06,bs}.}
      if the condition 

\medskip

\noindent $(\mathsf{H}) {}~~~~~~~~ $
     $ \ \ \,\_f(\U)= \_g(\U)\ \Longleftrightarrow \
      \{\,x\in X \mid f(x)=g(x)\} \in \U .$

\medskip      
\noindent holds for all $f,g:X\to X$.

\smallskip

So the condition $(\mathsf{H})$ has to be verified by all ultrafilters $\U_\xi$ in order that $\X$ be Hausdorff. We have thus proved the following
\begin{thm}\label{lpd}
  A  functional
      extension $\X$ of $X$ satisfies simultaneously the principles
     \emph{\ind\ }and \emph{\tran\ }if and
only if all the ultrafilters $\U_{\xi}$ on $X$ associated to the points  $\xi\in\X$ are Hausdorff. 
\qed
\end{thm}
\smallskip
We shall deal in the next section with the set theoretic strength of 
the combination of \ind\ with \tran. By now we simply recall that 
there are plenty of non-Hausdorff ultrafilters on $X$: \pes\ all  \emph{tensor products}
  $\V=\U\otimes\U$ of an ultrafilter $\U$ with itself contradict $(\mathsf{H})$, because  the ``projections" $p_1,p_2$ of the property \dir\  give $\_p_1(\V)=\U=\_p_2(\V)$.
  
   Thus we can easily state what we announced in the introduction, namely
\begin{thm}\label{lp123}
 No 
   extension satisfies at once the three Leibnizian principles
   \emph{\ind, \poss}, and \emph{\tran}.
\hfill $\Box$
\end{thm}

\section{Set theoretic problems in combining \ind\ with \tran.}\label{shaus}

As shown by Theorem \ref{lpd}, combining \ind\ with \tran\
requires  special ultrafilters, named \emph{Hausdorff} in 
Subsection \ref{trin}. Despite the apparent weakness of their defining 
property $(\mathsf{H})$, which is actually true whenever any of
the involved functions is i\emph{njective} (or \emph{constant}), not much is known about 
Hausdorff ultrafilters.

On  \emph{countable} sets, the property
$(\mathsf{H})$ is satisfied by \emph{selective}
ultrafilters as well as by \emph{products of pairwise non-isomorphic 
selective}
ultrafilters (see \cite{DF06}), but their existence in pure \zfc\ is still unproved.
However any hypothesis providing infinitely many 
non-isomorphic 
selective ultrafilters over $\N$, like the Continuum Hypothesis 
$\mathsf{CH}$ or Martin's Axiom $\mathsf{MA}$,
provides any \emph{countable} set with infinitely many non-isomorphic functional extensions  that 
satisfy \ind.

On \emph{uncountable} sets the situation is highly problematic: 
it is proved in \cite{DF06} that  Hausdorff ultrafilters on sets of size
greater than or equal to $\uf$ (the least size of an ultrafilter basis on $\N$)
cannot be \emph{regular}. All what is provable in 
$\mathsf{ZFC}$ about the size of $\uf$
is that $\aleph_{1} \le \uf \le 2^{\aleph_{0}}$ 
(see e.g. \cite{bl03}). 
In particular, the existence of functional extensions
satisfying \ind\ with uniform 
ultrafilters, even on $\R$, would  imply that of inner models 
with measurable 
 cardinals. (Such ultrafilters have been obtained by forcing only assuming
 much stronger hypotheses, see \cite{kt}).

Be it as it may, as far as we do not abide $\mathsf{ZFC}$ as our foundational
theory, \emph{the existence of
functional extensions without indiscernibles, although consistent, cannot be proved}.

\subsection{The choice of Gottfried Wilhelm von Leibniz}\label{lchoi}

  We have seen that (at least) one of the three Leibnizian principles that we 
  have investigated has to be left out. The most reasonable choice seems 
  to be that of dropping \ind. In fact, even if one neglects the set 
  theoretic problems that have been outlined above, one should pay 
  attention to what 
  Leibniz, whose logico-mathematical insight into philosophical questions cannot be overestimated, wrote about it.
  
   \begin{quote}

[...] cette supposition de deux indiscernables [...] paroist possible 
en termes abstraits, mais elle n'est point compatible avec l'ordre 
des choses [\ldots]

Quand je nie qu'il y ait [\ldots] deux corps indiscernables, je ne 
dis point qu'il soit impossible 
absolument d'en poser, mais que c'est une chose contraire a la 
sagesse divine [\ldots]

Les parties du temps ou du lieu [...] sont des choses ideales, ainsi
elles se rassemblent parfaitement comme deux unit\'es abstraites. 
Mais il n'est pas de m\^eme de deux Uns concrets [...] c'est \`a 
dire veritablement actuels.

Je ne dis pas que deux points de l'Espace sont un meme point, ny 
que deux instans du temps sont un meme instant comme il semble
qu'on m'impute [\ldots]\footnote{~
\ldots this supposition of two indiscernibles \ldots seems 
abstractly possible, but it is incompatible with the order of 
things\ldots\\
When I deny that there are  \ldots two indiscernible bodies, I do not 
say that [this existence] is absolutely impossible to assume, but 
that it is a thing contrary to divine wisdom \ldots \\
The parts of time or place \ldots are ideal things, so 
they perfectly resemble like
two abstract unities. But it is not so with two concrete Ones,\ldots
that is truly actual [things].\\
I don't say that two points of Space are one same point, neither 
that two instants of time are one same instant as it seems that one 
imputes to me \ldots}\
(\emph{Fifth letter to Clarke}, \cite{LC}, pp. 131-135)
   \end{quote}

   From these quotes, it appears that  Leibniz  himself considered the identity of indiscernibles 
   as a ``physical'' rather than a ``logical'' principle: it may be
   actually true, but its negation is non-contradictory in principle,
   so \emph{it could fail in some 
   possible world}. 
   Moreover only ``properties of the real world'' $M$
   are considered in  this principle: so it seems natural, and 
   not absurd, to assume that objects indiscernible by these ``real'' 
   properties may be separated by some abstract, ``ideal'' property 
   of $\hM$.
   
On this ground we finally decide  to call \emph{Leibnizian} a functional 
extension that satisfies both \textbf{Poss} and \textbf{Tran}, and so
necessarily not \textbf{Ind}. Topologically, this choice means that the $S$-topologies of the Leibnizian extensions $\hM$ are \emph{quasi-compact}, but only their restrictions to the ``standard' model $M$ are obviously Hausdorff (actually \emph{discrete}).  Thus the existence of plenty of 
\emph{Leibnizian extensions} is granted by the final results of Section 
\ref{stop}, without any need of supplementary set theoretic hypotheses.

\section{Combining \ind\ with \poss.}\label{inos}
By Theorems \ref{eltop} and \ref{lp123}, the principles \ind\ and \poss\ cannot hold simultaneously in a functional extension, where \tran\ always holds. So, if we want to justify our initial assertion that the three Leibnizian principles may be verified in pairs, we have to relax   the properties of the $\ast$-extension given in Definition \ref{text}, so as to allow the $S$-topology to be \emph{compac}t. 

%Taking into account that we are interested in extensions whose $S$-topology is Hausdorff, we know that the extended functions Have to be continuous, hence uniquely determined. whose

Of course, such a weakening should maintain most preservation properties of the functional extensions, although necessarily losing part of the \trans. In particular, by compactness, \emph{all} ultrafilters on $X$ have to be realized as associated to a \emph{unique} point $\xi\in\X$. It follows that the typical ``analytic" property of  the \emph{nonstandard} extension of functions, namely that 
\begin{center}
\emph{everywhere different functions have  everywhere different extensions, }
\end{center}
has to \emph{fail} in the presence of non-Hausdorff ultrafilters.

Fortunately, such a weakening may be realized still maintaining all the typical preservation properties of the \emph{continuous} extensions of functions, which are allowed to reach equality \emph{only at limit points}. 
In fact the crucial property of Subsection \ref{trin} \\
\noindent $(\mathsf{H}) {}~~~~~~~~ $
     $ \ \ \,\_f(\U)= \_g(\U)\ \Longleftrightarrow \
      \{\,x\in X \mid f(x)=g(x)\} \in \U .$\\
 holds for \emph{arbitrary} ultrafilters 
whenever at least one of the functions $f,g:X\to X$  either is \emph{injective} or has \emph{finite range}. So we  may maintain the conditions \comp\ and \dir, and weaken  \equ\ by postulating it only in those cases.

Therefore we define a \emph{weak functional extension} as follows:
\begin{defi}\label{wext}
A
\emph{weak functional extension} of $X$ is  a  superset $\X$ of $X$ where a
\emph{distinguished extension} $\f:\X \to \X$ is assigned to
each function $f:X \to X$, so as to satisfy the following conditions
%for all $f,g:X \to X$. 
%and all $\xi\in\X$:}
\begin{description}
\item[$(\mathsf{comp})$]  \textbf{Preservation of compositions:} $\! \ns(g\circ f)= \g\circ\!\f$ for all $f,g:X \to X$;
%\bigskip
\item[$(\mathsf{wequ})$] \textbf{Weak preservation of equalizers:}~$\ns(\chi_{fg})= \chi_{\f\g}\ $
 if at least one of the functions $f,g:X\to X$  either is injective or has finite range;
\item [$(\mathsf{dir})$] \textbf{Directness of the Puritz order:}
for all $\xi,\eta \in \X$ there
exists  $\zeta \in \X$ such that
\  $\xi,\eta\lesssim_P \zeta$.
\end{description}
\end{defi}
Also in this weaker context,  the  extension $\ns\chi_A$ of the characteristic function $\chi_A$ of a subset $A\incl X$ can be taken as the characteristic function of
the $\ast$-{\em extension} $\hA$ of
$A$ in $\X$, thus obtaining   again a 
\emph{Boolean isomorphism} between the field $\P(X)$ of all subsets of $X$ 
and a field $\C\ell(\X)$ of subsets of $\X$. So the latter can again generate the $S$-topology, with all the properties stated in  Section \ref{stop}; in particular the $\ast$-extended functions are continuous (and unique when \ind\ holds), and many important instances of \tran\ can be deduced.

A detailed study of the weak functional extensions lies outside the scope of this article. We simply recall that \ind\ and \poss\ together imply that the $S$-topology is compact, and that the map $\xi\mapsto \U_\xi$ establishes a \emph{biunique correspondence} between the points of the \emph{compact} extension $\X$ and the set of all ultrafilters over $X$.
The latter set can be naturally identified with the Stone-$\check{\rm C}$ech
compactification $\bg X$ of the discrete space $X$, with its usual topology having as basis
$\{ {\mathcal O}_{A} \mid A\in \P (X) \}$, where ${\mathcal O}_{A} $ is the
set of all ultrafilters containing A. (The embedding  $\, e:X \to \beta
X$ being given by the principal ultrafilters.) 

%$\xi\mapsto \U_\xi$ is biunique, and so it is a homeomorphism onto the Stone-$\check{\rm C}$ech
%compactification $\bg X$, if the latter is identified with the space of all ultrafilters on $X$. 

  In fact it turns out that the Stone-$\check{\rm C}$ech
compactification is essentially the unique weak functional extension of $X$ that satisfies both principles \ind\ and \poss, according to the following theorem, whose proof specializes that of Theorem 2.1 of \cite{DF05}.
 
\begin{thm}\label{indpos}
Let $\X$ be a weak functional extension of $X$ satisfying both principles \emph{\ind} and \emph{\poss}, endowed with the $S$-topology, and identify  the 
   Stone-$\check{\rm C}$ech
compactification $\bg X$ of the discrete space $X$ with the set of all ultrafilters on $X$. 
Then the map $\yg:\X \to\bg X$ defined by $$\yg(\xi)=\U_\xi=\{A\incl X\mid \xi\in\hA\}$$
establishes a homeomorphism  between $\X$ with the $S$-topology and  $\bg X$ with its compact topology. Moreover, for all $f:X\to X$, one has $$\_f\circ\yg=\yg\circ\f, $$
where $\_f$ is the unique continuous extension of $f$ to $\bg X$.
\end{thm}
\proof
~For all $x\in X$, $\U_{x}$ is the principal ultrafilter
generated by $x$, hence $\yg$ induces the canonical
embedding of  $X$ into $\beta X$. 

Moreover, for all $A\incl X$, the map $\yg$ induces a mapping of the clopen subset $\hA$
of $\,\X$ onto the basic
open set ${\mathcal O}_{A} $ of $\beta X$. Therefore $\yg$ is continuous and open with respect to the
$S$-topology of $\X$, \ie\
 $\yg $ is bicontinuous
between the $S$-topology of $\X$ and the compact topology of $\bg X$.

 On the other hand, the map $\yg$ is injective because  the
      $S$-topology is Hausdorff, and surjective because it is  quasi-compact.  So, being bicontinuous,   the map $\yg$ is  a homeomorphism.

    Finally, we have ~$\xi \in \hA \ \Leftrightarrow \
         \f(\xi) \in \st(f(A))$, for all $\, \xi \in \X$,  or equivalently
         $A\in {\U}_{\xi} \ \Leftrightarrow \ f(A) \in {\U}_{\f(\xi)}$.
         It follows that $\overline f \circ \yg= \yg \circ \f$,
         and the last assertion of the theorem is proved.
\qed

\medskip

\section{Final Remarks and Open Questions}\label{froq} 
\bigskip
Remark that we have proved that all Hausdorff (weak) functional extensions use  the
same ``function-extending mechanism'', namely that arising
from the Stone-$\check{\rm C}$ech
    compactification. Therefore, in the Hausdorff case, the choice of the
$\ast$-extensions of functions is forced by the unique topological compactification of $X$.

Also remark that the conditions $(\mathsf{equ})$ and $(\mathsf{dir})$ 
are independent, 
even when \ind\ holds, provided that  Hausdorff ultrafilters exist. 
Call \emph{invariant} a subspace $Y$ of $\,\X$ (or of the Stone-$\check{\rm C}$ech compactification $\bg X$ of X)
%    if it isclosed under *functions, i.e. 
   if  
$$\f(\xi)\in Y\ \mbox{(respectively $\_f(\U_\xi)\in Y$) for all}\ f:X\to X\ 
\mbox{and all} \
\xi\in Y.$$

It is easily seen that any invariant subspace $Y$ of a (weak) functional extension $\,\X$ is itself a (weak)
functional extension of $X$, by taking the restrictions to $Y$ of  the functions $\f$ for all $f:X\to X$; and clearly the corresponding ultrafilters $\U_\xi, \,\xi\in Y$ constitute an
invariant  subspace of $\beta X$.

Now 
there are invariant 
subspaces of  $\bg X$
where $(\mathsf{equ})$ holds whereas 
$\mathsf{(\mathsf{dir})}$ fails and \emph{vice versa}, as well as  
invariant 
subspaces where both fail or hold. 
In fact, for $\U\in\bg X$ let
$Y_{\U}=\{\_f(\U)\mid f:X\to X\,\}$ be the invariant subspace 
generated by $\U$. Clearly $Y_{\U}$ is directed, so $(\mathsf{dir})$ holds in $Y_{\U}$ for 
all ultrafilters $\U$, whereas $(\mathsf{equ})$ holds if and only if $\,\U$ is 
Hausdorff.
On the other hand, let $\U$ and $\V$ be ultrafilters such 
that  neither 
of them belongs to the invariant subspace generated by the other one: 
then $Y_{\U}\cup Y_{\V}$ is an invariant subspace where $(\mathsf{dir})$ fails, while $(\mathsf{equ})$ holds if and only if both $\U$ and $\V$ are Hausdorff.

\subsection{Some open questions}\label{sopq}

We conclude the paper with three open questions that involve special
ultrafilters, and so should
 be of independent set theoretic interest.
\begin{enumerate}
    \item Is the existence of functional extensions of $\N$ \emph{without 
indiscernibles}
provable in $\mathsf{ZFC}$, or at least 
derivable from
set-theoretic hypotheses weaker than those providing \emph{selective 
ultrafilters}? \Pes\ from ${\mathfrak x = \mathfrak c}$, where  
$\mathfrak x$ is a
 cardinal 
invariant of the continuum not dominated by \textbf{cov}$(\B)$ (see \cite{bl03})?
  \item  Is it consistent with $\mathsf{ZFC}$ that there exist 
    nonstandard real lines 
$\hR$ \emph{without indiscernibles} where all points correspond to \emph{uniform} ultrafilters?
 \item  Is the existence of \emph{countably compact functional extensions}
consistent with $\mathsf{ZFC}$? (These extensions would be of great 
interest, because they would verify \ind, 
\tran, and the weakened version of \poss\ that considers only 
\emph{sequences} of properties.)
\end{enumerate}

\section*{Acknowledgment}
  \noindent {{The author is grateful to  Vieri Benci, Mauro Di Nasso and Massimo 
Mugnai for useful discussions
and suggestions. The author thanks the anonymous referees for several important remarks and suggestions.}}

%% The Appendices part is started with the command \appendix;
%% appendix sections are then done as normal sections
%% \appendix

%% \section{}
%% \label{}

%% References
%%
%% Following citation commands can be used in the body text:
%% Usage of \cite is as follows:
%%   \cite{key}          ==>>  [#]
%%   \cite[chap. 2]{key} ==>>  [#, chap. 2]
%%   \citet{key}         ==>>  Author [#]

%% References with bibTeX database:

% \bibliographystyle{model1-num-names}
% \bibliography{<your-bib-database>}

\begin{thebibliography}{99}

%% \bibitem must have the following form:
%%   \bibitem{key}...
%%

\bibitem{nsa} {\sc L.O.~Arkeryd, N.J.~Cutland, C.W.~Henson} (eds.) -
\emph{Nonstandard Analysis - Theory and Applications}. NATO ASI Series \textbf{C 493}, 
 Kluwer A.P., Dordrecht 1997.
  \bibitem{bs} \textsc{T.~Bartoszynski, S.~Shelah} - 
On the density of Hausdorff ultrafilters, in
 \emph{Logic Colloquium 2004}, L. N. Logic \textbf{29}, A. S. L., 
 Chicago 2008, 18--32.\
\bibitem{BDF02} {\sc V.~Benci, M.~Di~Nasso, M.~Forti} -
Hausdorff Nonstandard Extensions, \emph{Bol.~Soc.~Parana.~Mat.}
(3)\textbf{20} (2002), 9--20.
\bibitem{BDF06} \textsc{V.~Benci, M.~Di~Nasso, M.~Forti} -
The Eightfold Path 
to Nonstandard Analysis, in \emph{Nonstandard Methods and 
Applications in Mathematics} (N.J.~Cutland, M.~Di~Nasso, D.A.~Ross, 
eds.), L.N. in Logic \textbf{143}, A.S.L. 2006, 3--44.
\bibitem{bl03} \textsc{A.~Blass} - Combinatorial cardinal characteristics of the 
continuum, in  \emph{Handbook of Set Theory} (M. Foreman and  
A. Kanamori, eds.), Springer V., Dordrecht 2010, 395--490.
\bibitem{CK} \textsc{C.C.~Chang, H.J.~Keisler} -
\textit{Model Theory} (3rd edition).
North-Holland, Amsterdam 1990.
\bibitem{dt} \textsc{M.~Daguenet-Teissier} - Ultrafiltres \`a la
fa\c{c}on de Ramsey,
\textit{Trans. Amer. Math. Soc.}, \textbf{250} (1979), 91--120.
\bibitem{DF05} 
\textsc{M.~Di~Nasso, M.~Forti} - 
Topological and nonstandard extensions,
\textit{Monatsh. f. Math.} {\bf 144} (2005), 89--112.
\bibitem{DF06} 
{\sc M.~Di~Nasso, M.~Forti} -
Hausdorff Ultrafilters, \emph{Proc. Amer. Math. Soc.} \textbf{134}
(2006), 1809--18.
\bibitem{Eng} \textsc{R.~Engelking} - \emph{General Topology}. Polish
S.~P., Warszawa 1977.
\bibitem{F12} \textsc{M. Forti} - A simple algebraic characterization of nonstandard extensions, 
\emph{Proc.~Amer.~Math.~Soc.} {\bf 140} (2012), 2903--2912. 
\bibitem{kt} \textsc{A.~Kanamori, A.D.~Taylor} - {Separating
ultrafilters on uncountable
cardinals},
\emph{Israel J. Math.} \textbf{47} (1984), 131--138.
\bibitem{Ke76} \textsc{H.J.~Keisler} -
\textit{Foundations of Infinitesimal Calculus}.
Prindle, Weber and Schmidt, Boston 1976.
\bibitem{LC} {\sc G.~W.~Leibniz} - \emph{Correspondance Leibniz-Clarke, 
pr\'esent\'ee d'apr\`es les manuscrits originaux des biblioth\`eques 
de Hanovre et de Londres}, par Andr\'e Robinet, PUF, Paris 1957.
\bibitem{ng} {\sc S.-A.~Ng, H.~Render} - The Puritz order and its relationship to the Rudin-Keisler order, in \emph{Reuniting the antipodes - Constructive and nonstandard views of the continuum} (P.~Schuster, U.~Berger, H.~Oswald, eds.), Kluwer, Dordrecht 2001, 157--166.
         \bibitem{pu} {\sc C.~Puritz} - Ultrafilters and standard functions in nonstandard arithmetic, \emph{Proc.~London Math.~Soc} (3) \textbf{22} (1971), 705--733                         
\bibitem{rob66} \textsc{A.~Robinson} - \emph{Non-standard Analysis}. North 
Holland, Amsterdam 1966.
	  
% 
% 	   \vspace{-11mm}
%[L--E] LEIBNIZ--EDITION - http://www.leibniz-edition.de/


 \end{thebibliography}

%% Authors are advised to submit their bibtex database files. They are
%% requested to list a bibtex style file in the manuscript if they do
%% not want to use model1-num-names.bst.

%% References without bibTeX database:

\end{document}